\documentclass[12pt,reqno,twoside,letterpaper]{amsart}

\usepackage{hyperref}

\newcommand{\C}{\mathfrak{C}}
\newcommand{\N}{\mathbb{N}}

\newcommand{\Mod}[1]{\ (\mathrm{mod}\ #1)}

\begin{document}

\title[On the cluster structures in Collatz level sets]{On the cluster structures in Collatz level sets}
\author{Markus Sigg}
\address{Freiburg, Germany}
\email{mail@markussigg.de}
\date{December 20, 2020}

\begin{abstract}
  The cluster structures that can be observed in the first few level sets of
  the Collatz tree are maintained through all its levels, provided that the
  \emph{orbit steadiness}
  \begin{equation*}
    \prod_{\substack{k \in R(n)\\ k \equiv 4 \Mod 6}} \frac{k-1}k
  \end{equation*}
  of the elements $n$ of the Collatz tree is suitably bounded from below, where
  $R(n)$ denotes the Collatz orbit of $n$.
\end{abstract}

\maketitle

{
}

\section{The question}

Let $\N$ be the set of positive integers. By $c:\N \longrightarrow \N$ we
denote the \emph{Collatz function} (see \cite{lagarias}), defined by
\begin{equation*}
  c(n) := \begin{cases}
    \displaystyle \phantom{3n}\frac n2 & \text{if $n$ is even,}\\[1em]
    3n+1 & \text{if $n$ is odd.}
  \end{cases}
\end{equation*}
We are interested in the \emph{level sets}
$L_\nu := \{ n \in \N : c^\nu(n) = 1\ \wedge\ c^i(n) \neq 1\ \text{for}\ i < \nu\}$ for
$\nu \in \N_0 := \N \cup \{0\}$, whose elements are listed in \cite{noe}. It is
easy to see that
\begin{equation}\label{lnu}
  L_{\nu+1} =
  \left\{2n : n \in L_\nu\right\}
  \uplus
  \left\{\frac{n-1}3 : 4 < n \in L_\nu, n \equiv 4 \Mod 6\right\}.
\end{equation}
The level sets for the first some dozens values of $\nu$ exhibit a cluster
structure. For example, the $72$ elements of $L_{20}$ are given by the seven
clusters

\begin{center}
  $\{18, 19\},$\\
  $\{112, 116, 117, 120, 122\},$\\
  $\{704, 720, 724, 725, 736, 738, 739, 744, 746, 753, 802, 803, 804, 805, 806\},$\\
  $\{4352, \left<20\right>, 4849\},$\\
  $\{24576, \left<17\right>, 29126\},$\\
  $\{163840, 172032, 174080, 174592, 174720, 174752, 174760, 174762\},$\\
  $\{1048576\}$,
\end{center}

where $\left<x\right>$ is short for a list of $x$ intermediate elements. Of
course the largest number $1048576$ equals $2^{20}$, and (\ref{lnu}) makes
plausible that the elements of a cluster have about six times the size of the
elements of the previous cluster: Each cluster $C$ of $L_\nu$ gives a cluster
$2C$ of $L_{\nu+1}$, and the second set in (\ref{lnu}) contributes to the
cluster below that one. An additional cluster appears at the bottom end in
$L_{\nu+1}$ whenever there exists an $n \equiv 4 \Mod 6$ in the lowest cluster
of $L_\nu$. The question suggests itself whether the pattern extends to all
level sets or whether, for large $\nu$, the clusters eventually dissolve so
much that they overlap.

\section{The answer (under a provision)}

A closer investigation of the issue shows that the product of quotients
$\frac{k-1}k$ for orbit elements $k \equiv 4 \Mod 6$ controls the evolution of
the cluster shapes. We introduce the \emph{orbit steadiness} function
$\sigma:\N \longrightarrow [0,1]$ by
\begin{equation*}
  \sigma(n) := \prod_{\substack{k \in R(n)\\ k \equiv 4 \Mod 6}} \frac{k-1}k
  \quad \text{for $n \in \N$},
\end{equation*}
where $R(n) := \{c^i(n): i \in \N_0\}$ is the orbit set of $n$, and with
$\C := \bigcup\limits_{\nu \in \N_0} L_\nu$
\begin{equation*}
  \sigma_0 := \inf_{n \in \C} \sigma(n).
\end{equation*}
For $\nu \in \N_0$, by
\begin{equation*}
  S_{\nu,\kappa} := \left[ \sigma_0 \, \frac{2^\nu}{6^\kappa}, \frac{2^\nu}{6^\kappa} \right]
  \quad \text{for $\kappa \in \N_0$}
\end{equation*}
we define ``slots'' for the clusters of the level set $L_\nu$. Then, given
$\nu \in \N_0$ and $n \in L_\nu$, setting
$I_\alpha := \{ i \in \{1, \dots, \nu\} : c^{i-1}(n) \equiv \alpha \Mod 2 \}$
for $\alpha \in \{0,1\}$ and $\kappa := |I_1|$ yields
\begin{equation*}
  n = \prod_{i=1}^\nu \frac{c^{i-1}(n)}{c^i(n)}
  = \prod_{i \in I_0} \frac{c^{i-1}(n)}{c^i(n)} \cdot \prod_{i \in I_1} \frac{c^{i-1}(n)}{c^i(n)}
  = 2^{|I_0|} \prod_{i \in I_1} \frac{\frac{c^i(n)-1}3}{c^i(n)}
  = \frac{2^{\nu - \kappa}}{3^\kappa} \, \sigma(n)
  = \frac{2^\nu}{6^\kappa} \, \sigma(n),
\end{equation*}
which gives $n \le \frac{2^\nu}{6^\kappa}$ and
$n \ge \sigma_0 \, \frac{2^\nu}{6^\kappa}$, hence $n \in S_{\nu,\kappa}$. This
proves
\begin{equation*}
  L_\nu \subset \bigcup_{\kappa \in \N_0} S_{\nu,\kappa}
  \quad \text{for all $\nu \in \N_0$.}
\end{equation*}
For $\nu \in \N _0$, the slots $S_{\nu,\kappa}$ for $\kappa \in N_0$ are
pairwise disjoint if $\sigma_0 > \frac16$. Numerical evidence indicates
$\sigma_0 \approx 0.5152$. The slightly weaker assumption $\sigma_0 > \frac12$
gives $\max S_{\nu,\kappa} < \frac13 \min S_{\nu,\kappa -1}$ for all
$\kappa \in \N$, i.\,e.\ a clear separation of the slots and the persistence of
the cluster pattern in all level sets of the Collatz tree.

To get a trustworthy statement, a proven lower bound for $\sigma_0$ would of
course be preferable.

\section{A remark}

We admitted only elements of the Collatz tree in the definition of $\sigma_0$.
This does not make a difference if the Collatz conjecture is true, because then
$\C = \N$. However, at this point it cannot be excluded that, in case of the
falsehood of the Collatz conjecture, the orbit steadiness of some
$n \in \N \setminus \C$, i.\,e.\ some $n$ with a non-trivial cyclic or a
diverging orbit, might be smaller than $\sigma_0$. This is why, in an abundance
of caution, and because it is sufficing for the application, we decided for
$\C$ instead of $\N$.

\end{document}